\documentclass[a4paper, 11pt]{amsart}

\usepackage[T1]{fontenc}

\usepackage[abbrev]{amsrefs}

\theoremstyle{plain}
\newtheorem{proposition}{Proposition}[section]
\newtheorem{theorem}{Theorem}
\newtheorem{lemma}[proposition]{Lemma}

\theoremstyle{definition}

\newcommand{\N}{\mathbb{N}}
\newcommand{\R}{\mathbf{R}}

\newcommand{\abs}[1]{\mathopen\lvert#1\mathclose\rvert}
\newcommand{\st}{\,:\,}
\newcommand{\muleb}[1]{\mathcal{L}^{#1}}
\newcommand{\muhaus}[1]{\mathcal{H}^{#1}}

\DeclareMathOperator{\Div}{div}
\DeclareMathOperator{\supp}{supp}

\title[Extremal functions in Poincar\'e--Sobolev inequalities]{Extremal functions in Poincar\'e--Sobolev inequalities for functions of bounded variation}
\author{Vincent Bouchez}

\address{
 Universit{\'e} catholique de Louvain\\
 D\'epartement de Math{\'e}matique\\
 Chemin du cyclotron 2\\
1348 Louvain-la-Neuve\\
Belgium}
\email{Vincent.Bouchez@uclouvain.be}

\author{Jean Van Schaftingen}
\address{
Universit{\'e} catholique de Louvain\\
 D\'epartement de Math{\'e}matique\\
 Chemin du cyclotron 2\\
1348 Louvain-la-Neuve\\
Belgium}
\email{Jean.VanSchaftingen@uclouvain.be}
\subjclass[2000]{Primary: 46E35; Secondary: 26B30, 26D15, 35J62}
\keywords{Poincar\'e--Sobolev inequality, sharp constant, optimal constant, extremal function, function of bounded variation, concentration--compactness, scalar curvature, compact manifold, mean curvature, Gauss-Bonnet formula}

\dedicatory{Dedicated to Jean--Pierre Gossez, on the occasion of his 65th birthday}

\begin{document}

\begin{abstract}
If \(\Omega \subset \R^n\) is a smooth bounded domain and \(q \in (0, \frac{n}{n-1})\) we consider the Poincar\'e--Sobolev inequality
\[
  c \Bigl(\int_{\Omega} \abs{u}^\frac{n}{n-1}\Bigr)^{1-\frac{1}{n}} \le \int_{\Omega} \abs{Du},
\]
for every \(u \in \mathrm{BV}(\Omega)\) such that \(\int_{\Omega} \abs{u}^{q-1} u = 0\). We show that the sharp constant is achieved.
We also consider the same inequality on an \(n\)--dimensional compact Riemannian manifold \(M\). When \(n \ge 3\) and the scalar curvature is positive at some point, then the sharp constant is achieved. In the case \(n \ge 2\), we need the maximal scalar curvature to satisfy some strict inequality.
\end{abstract}

\maketitle

\section{Introduction}

If \(\Omega \subset \R^n\) is smooth and has finite measure and if \(p \in (1, n)\), there exists \(c > 0\) such that for every \(u\) in the Sobolev space \(\mathrm{W}^{1, p}(\Omega)\) with \(\int_{\Omega} u = 0\), 
\begin{equation}
\label{ineqSobolevAverage}
 c \Bigl(\int_{\Omega} \abs{u}^\frac{np}{n-p}\Bigr)^{1-\frac{p}{n}} \le \int_{\Omega} \abs{\nabla u}^p.
\end{equation}
This inequality follows from the classical Sobolev inequality 
\[
  \Bigl( \int_{\Omega} \abs{u}^\frac{np}{n-p} \Bigr)^{1-\frac{p}{n}} \le C\Bigl(\int_{\Omega} \abs{\nabla u}^p+\abs{u}^p\Bigr)
\]
and a standard compactness argument (see for example E.\thinspace Giusti \cite{Giusti2003}*{\S 3.6}).

We are interested in whether the sharp constant in \eqref{ineqSobolevAverage}, that is the value
\[
\inf \Bigl\{\int_{\Omega} \abs{\nabla u}^p \st u \in \mathrm{W}^{1,p}(\Omega), \int_{\Omega} \abs{u}^{\frac{np}{n-p}} = 1 \text{ and } \int_{\Omega} u=0 \Bigr\},  
\]
is achieved. Since the embedding of \(\mathrm{W}^{1,p}(\Omega)\) in \(\mathrm{L}^{\frac{np}{n-p}}(\Omega)\) is not compact, the solution to this problem is not immediate. 

In the case where \(\Omega = \R^n\) and the condition \(\int_{\Omega} u = 0\) is dropped, this was solved by T.\thinspace Aubin~\cite{Aubin1975} and G.\thinspace Talenti \cite{Talenti1976}.
When the condition \(\int_{\Omega} u = 0\) is replaced by \(u=0\) on \(\partial \Omega\), it is known that the constant is not achieved. However, if \(u\) is only required to vanish on a part \(\Gamma\) of the boundary and \(\Gamma\) has some good geometric properties, P.-L.\thinspace Lions, F. Pacella and M. Tricarico have showed that the corresponding sharp constant is achieved for every \(p \in (1, \Bar{p})\) where \(\Bar{p} \in (1, n]\) depends on \(\Omega\) and \(\Gamma\) \cite{LionsPacellaTricarico1988}. Returning to our problem  P.\thinspace Gir\~ao and T.\thinspace Weth \cite{GiraoWeth2006} have showed that the sharp constant is achieved for \(p=2\). A.\thinspace V.\thinspace Demyanov and A.\thinspace I.\thinspace Nazarov~\cites{DemyanovNazarov2005, Nazarov2008} have proved that there exists \(\delta>0\) depending on \(\Omega\) such that the sharp constant is achieved for \(p \in (1, \frac{n+1}{2}+\delta)\). M.\thinspace Leckband \cite{Leckband2009} has given an alternative proof of this statement for a ball. 

\bigbreak

We are interested in the same question when \( p=1 \).  
The counterpart of the Sobolev space \(\mathrm{W}^{1, p}(\Omega)\) in this case is the space of functions of bounded variation \(\mathrm{BV}(\Omega)\), and the inequality~\eqref{ineqSobolevAverage} becomes
\begin{equation}
 c \Bigl(\int_{\Omega} \abs{u}^\frac{n}{n-1}\Bigr)^{1-\frac{1}{n}} \le \int_{\Omega} \abs{D u},
\end{equation}
where now \(\abs{Du}\) is a measure. The sharp constant is then
\[
  c^1_\Omega=\inf \Bigl\{\int_{\Omega} \abs{D u} \st u \in \mathrm{BV}(\Omega), \int_{\Omega} \abs{u}^{\frac{n}{n-1}} = 1 \text{ and } \int_{\Omega} u=0 \Bigr\}.
\]
When \(\Omega\) is a ball, A.\thinspace Cianchi~\cite{Cianchi1989} has showed that the sharp constant is achieved. In the general case, Zhu M. \cite{Zhu2004CCM}*{Theorem 1.3} has showed that if one restricts the inequality to functions in \(\mathrm{BV}(\Omega)\) that take two values, the sharp constant is achieved.

Our first result is 

\begin{theorem}
\label{theoremDomainAverage}
Let \(n \ge 2\). If \(\Omega \subset \R^n\) is a bounded domain with \(\mathcal{C}^2\) boundary, then \(c^1_\Omega\) is achieved.
\end{theorem}

When \(n=2\), this answers a question mentioned by H.\thinspace Brezis and J.\thinspace Van Schaftingen \cite{BrezisVanSchaftingen2008}*{problem 3}.

\bigbreak

Instead of considering the inequality~\eqref{ineqSobolevAverage} under the constraint \(\int_{\Omega} u = 0\), one can drop the condition and take the infimum over functions that only differ by a constant:
\begin{equation}
\label{ineqSobolevMinimal}
 c \inf_{\lambda \in \R} \Bigl(\int_{\Omega} \abs{u-\lambda}^\frac{np}{n-p}\Bigr)^{1-\frac{p}{n}} \le \int_{\Omega} \abs{\nabla u}^p;
\end{equation}
this is equivalent to \eqref{ineqSobolevAverage} under the constraint \(\int_{\Omega} \abs{u}^{\frac{np}{n-p}-2} u=0\). In this setting, A.\thinspace V. Demyanov and A.\thinspace I.\thinspace Nazarov~\cite{DemyanovNazarov2005}*{theorem~7.4} have proved that when
\begin{multline*}
  1 < p < \max\Bigl(\frac{3n+1-\sqrt{5n^2+2n+1}}{2},\\
 \frac{n^2+3n+1+\sqrt{n^4+6n^3-n^2-2n+1}}{2(3n+2)}\Bigr),
\end{multline*}
the optimal constant
\[
\inf \Bigl\{\int_{\Omega} \abs{\nabla u}^p \st u \in \mathrm{W}^{1,p}(\Omega), \int_{\Omega} \abs{u}^{\frac{np}{n-p}} = 1 \text{ and } \int_{\Omega} \abs{u}^{\frac{np}{n-p}-2}u=0 \Bigr\}
\]
is achieved. 

We consider the corresponding problem of determining whether
\[
c^\frac{1}{n-1}_\Omega=\inf \Bigl\{\int_{\Omega} \abs{D u} \st u \in \mathrm{BV}(\Omega), \int_{\Omega} \abs{u}^{\frac{n}{n-1}} = 1 \text{ and } \int_{\Omega} \abs{u}^{\frac{1}{n-1}-1}u=0 \Bigr\}  
\]
is achieved. 

\begin{theorem}
\label{theoremDomainMinimal}
Let \(n \ge 2\). If \(\Omega \subset \R^n\) is a bounded domain with \(\mathcal{C}^2\) boundary, then \(c^{\frac{1}{n-1}}_\Omega\) is achieved.
\end{theorem}

More generally, we can consider the quantity
\[
  c^q_\Omega = \inf \Bigl\{\int_{\Omega} \abs{D u} \st u \in \mathrm{BV}(\Omega), \int_{\Omega} \abs{u}^{\frac{n}{n-1}} = 1 \text{ and } \int_{\Omega} \abs{u}^{q-1}u=0 \Bigr\}
\]
for \(q \in (1, \frac{n}{n-1})\).
If \(\Omega \subset \R^n\) is a bounded domain with \(\mathcal{C}^1\) boundary, \(c^q_\Omega > 0\).
In this setting theorems~\ref{theoremDomainAverage} and \ref{theoremDomainMinimal} are particular cases of 

\begin{theorem}
\label{theoremDomainGeneral}
Let \(n \ge 2\). If \(\Omega \subset \R^n\) is a bounded domain with \(\mathcal{C}^2\) boundary, then for every \(q \in (0, \frac{n}{n-1})\), \(c^{q}_\Omega\) is achieved.
\end{theorem}

\bigbreak

The inequality \eqref{ineqSobolevAverage} is also valid on a compact manifold without boundary \(M\). In this setting, 
Zhu M.\thinspace  \cites{Zhu2004PJM,Zhu2004NEEDS} has showed that the sharp constant is achieved when \(p \in (1, (1+\sqrt{1+8n})/4)\) on the sphere. A.\thinspace V.\thinspace Demyanov and A.\thinspace I.\thinspace Nazarov have showed that if there exist a point of \(M\) at which the scalar curvature is positive, then then there exists \(\delta>0\) such that the corresponding sharp constant is achieved for \(n \ge 3\) and \(p \in (1, \frac{n+2}{3}+\delta)\)~\cite{DemyanovNazarov2005}*{theorem 5.1}. For the inequality \eqref{ineqSobolevMinimal}, 
A.\thinspace V.\thinspace Demyanov and A.\thinspace I.\thinspace Nazarov~\cite{DemyanovNazarov2005}*{theorem 6.1} have proved that if the scalar curvature is positive at some point and
\begin{multline*}
 1 < p < \max\Bigl(2n+1-\sqrt{3n^2+2n+1},\\
 \frac{n^2+6n+2+\sqrt{n^4+12n^3-8n+4}}{2(5n+4)}\Bigr),
\end{multline*}
then the sharp constant is achieved\footnote{In some cases, the condition on the scalar curvature is reversed. Considers the quantity
\[
  \sup \Bigl\{ \Bigl(\int_{M} \abs{u}^{\frac{2n}{n-2}}\Bigr)^{1-\frac{2}{n}} - K_{2, n}\int_{M} \abs{\nabla u}^2 \st  \int_M \abs{u} = 1\Bigr\},
\]
where \(K_{2, n}=\frac{1}{n (n-2) \pi} (\frac{\Gamma(n)}{\Gamma(n/2)})^\frac{2}{n}\) is the optimal Sobolev constant on \(\R^n\) \citelist{\cite{Aubin1975}\cite{Talenti1976}}
; if \(n \ge 4\) the supremum is finite and achieved if \(M\) has \emph{negative} scalar curvature \cite{Hebey}*{theorem 1} and whereas it is not finite if the scalar curvature is positive somewhere \cite{DruetHebeyVaugon}*{theorem 1}.}. Moreover, they have proved that for \(p \ge \frac{n+1}{2}\) the sharp constant is not achieved on the \(n\)--dimensional sphere\footnote{This is in contradiction with a result of Zhu M.~\cite{Zhu2004NEEDS}*{theorem 3.1}.}.

For a compact \(\mathcal{C}^1\) Riemannian manifold \(M\) of dimension \(n \ge 2\) we consider whether the quantity
\[
  c^q_M = \inf \Bigl\{\int_{\Omega} \abs{\nabla u}^p \st u \in \mathrm{BV}(\Omega), \int_{\Omega} \abs{u}^{\frac{n}{n-1}} = 1 \text{ and } \int_{\Omega} \abs{u}^{q-1}u=0 \Bigr\}
\]
is achieved, with \(q \in (0, \frac{n}{n-1})\).
In the case where \(n \ge 3\) and the manifold has somewhere positive scalar curvature, one has the counterpart of theorem~\ref{theoremDomainGeneral}

\begin{theorem}
\label{theoremnManifolds}
Let \(n \ge 3\) and \(M\) be an \(n\)--dimensional compact Riemannian \(\mathcal{C}^2\) manifold. If there exists \(a \in M\) such that the scalar curvature \(S_a\) at \(a\) is positive, then for every \(q \in (0, \frac{n}{n-1})\), \(c^q_M\) is achieved.
\end{theorem}

In dimension \(2\), the same method only yields

\begin{theorem}
\label{theorem2ManifoldsSubcritical}
Let \(M\) be a \(2\)--dimensional compact Riemannian \(\mathcal{C}^2\) manifold. If there exists \(a \in M\) such that the scalar curvature \(S_a\) at \(a\) is positive, then for every \(q \in (0, 1)\), \(c^q_M\) is achieved.
\end{theorem}

If we strengthen the condition on the curvature we obtain

\begin{theorem}
\label{theorem2ManifoldsCritical}
Let \(M\) be a \(2\)--dimensional compact Riemannian \(\mathcal{C}^2\) manifold. If there exists \(a \in M\) such that the scalar curvature \(S_a\) at \(a\) satisfies
\[
 S_a > \frac{8\pi}{\muhaus{2}(M)},
\]
then for every \(q \in (0, 2)\), \(c^q_M\) is achieved.
\end{theorem}

Here $\muhaus{2}(M)$ denotes the two-dimensional Hausdorff measure of the manifold $M$.

In particular, theorem~\ref{theorem2ManifoldsCritical} allows to solve completely the case of surfaces of Euler--Poincar\'e characteristic \(2\) of nonconstant gaussian curvature.

\begin{theorem}
\label{theoremEulerPoincare}
Let \(M\) be a \(2\)--dimensional compact \(\mathcal{C}^2\) Riemannian manifold with nonconstant scalar curvature. If \(\chi(M)=2\), then for every \(q \in (0, 2)\), \(c^q_M\) is achieved.
\end{theorem}

While the sphere does not satisfy the hypotheses of the previous theorem, we have

\begin{theorem}
\label{theoremSphere}
For every \(q \in (0, 2)\), \(c^q_{\mathbf{S}^2}\) is achieved.
\end{theorem}

\section{Preliminaries}

Recall that for \(\Omega \subseteq \R^n\) open, \(\mathrm{BV}(\Omega)\) denotes the space of functions \(u \in \mathrm{L}^1(\Omega)\) such that 
\[
   \sup \Bigl\{ \int_{\Omega} u \Div \varphi \st \varphi \in \mathcal{C}^1_c(\R^n; \R^n) \text{ and } \abs{\varphi} \le 1\Bigr\} < \infty.
\]
If \(u \in \mathrm{BV}(\Omega)\), then there exists a vector measure \(Du\) such that for every \(\varphi \in \mathcal{C}^1_c(\R^n; \R^n)\),
\[
  \int_{\Omega} u \Div \varphi= -\int_{\Omega} \varphi \cdot Du.
\]
In particular, one can consider the variation \(\abs{Du}\) of \(Du\) which is a bounded measure on  \(\Omega\).

The optimal Sobolev inequality of H.\thinspace Federer and W.\thinspace H.\thinspace Fleming \cite{FedererFleming1960} states that for every \(u \in \mathrm{BV}(\R^n)\), 
\begin{equation}
\label{SobolevRn}
  \frac{\pi^{\frac{1}{2}} n}{\Gamma(\frac{n}{2}+1)^\frac{1}{n}} \Bigl( \int_{\R^n} \abs{u}^{\frac{n}{n-1}} \Bigr)^{1-\frac{1}{n}} \le  \int_{\R^n} \abs{Du}.
\end{equation}
The proof also shows that the constant if optimal and that it is achieved by multiples of characteristic functions of balls (see also \cite{CorderoNazaretVillani2004}).
If \(\R^n_+\) denotes the \(n\)--dimensional half-space, one deduces from \eqref{SobolevRn} by a reflexion argument that for every \(u \in \mathrm{BV}(\R^n_+)\), one has 
\[
 \frac{\pi^{\frac{1}{2}} n}{2^\frac{1}{n}\Gamma(\frac{n}{2}+1)^\frac{1}{n}} \Bigl( \int_{\R^n_+} \abs{u}^{\frac{n}{n-1}} \Bigr)^{1-\frac{1}{n}} \le  \int_{\R^n_+} \abs{Du}.
\]
One can show that the constant is achieved by characteristic functions of intersections of balls centered on the boundary of \(\R^n_+\) with \(\R^n_+\) itself.

A consequence that we shall use is 
\begin{lemma}
\label{lemmaInequality}
Let \(\Omega \subset \R^n\) be a domain with a \(\mathcal{C}^1\) boundary. 
For every \(a \in \partial \Omega\) and \(\varepsilon > 0\), there exists \(\delta > 0\) such that if \(u \in \mathrm{BV}(\Omega)\) and \(\supp u \subset B(a, \delta)\), then
\[
  \Bigl(\frac{c^*_n}{2^\frac{1}{n}}-\varepsilon\Bigr)\Bigl(\int_{\Omega} \abs{u}^{\frac{n}{n-1}}\Bigr)^{1-\frac{1}{n}}
\le \int_{\Omega} \abs{Du}.
\]
\end{lemma}

\section{Extremal functions on bounded domains}

\subsection{Existence by concentration-compactness}
A first ingredient in our proof of theorem~\ref{theoremDomainGeneral} is 

\begin{proposition}
\label{propositionSufficientDomain}
Let \(n \ge 2\) and \(\Omega \subset \R^n\) be a bounded domain with \(\mathcal{C}^2\) boundary. If
\[
  c^q_\Omega < \frac{c^*_n}{2^\frac{1}{n}},
\]
then \(c^q_\Omega\) is achieved.
\end{proposition}

In the case of the sharp constants for embeddings of \(\mathrm{W}^{1, p}(\Omega)\), with \(1 < p < n\), 
the counterpart has been proved has been proved by A.\thinspace V.\thinspace Demyanov and A.\thinspace I.\thinspace Nazarov~\cite{DemyanovNazarov2005}*{proposition 7.1}. An alternative argument has been provided by S.\thinspace De\thinspace Valeriola and M.\thinspace Willem \cite{DeValeriolaWillem2009}*{Theorem 4.1}.

Our main tool shall be
\begin{proposition}
\label{propositionCCLemma}
Let \((u_m)_{m \in \N}\) in \(\mathrm{BV}(\Omega)\) converge weakly to some \(u \in BV (\Omega)\).
Assume that there exist two bounded measures \(\mu\) and \(\nu\) on \(\Bar{\Omega}\) such that 
\((\abs{u_m}^{\frac{n}{n-1}})_{m \in \N}\) and \((\abs{Du_m})_{m \in \N}\) converge weakly in the sense of measures to \(\mu\) and \(\nu\) respectively.
Then there exists some at most countable set \(J\), distinct points \(x_j \in \Bar{\Omega}\) and real numbers \(\nu_j > 0\) with \(j \in J\) such that 
\begin{align*}
\nu&=\vert u\vert^{\frac{n}{n-1}}+\sum_{j \in J}\nu_j
\delta_{x_j},\\
\mu&\ge \abs{Du} + \frac{c^*_n}{2^\frac{1}{n}}\sum_{j\in
J}\nu_j^{1-\frac{1}{n}}\delta_{x_j}.
\end{align*}
\end{proposition}

This result is a variant of the corresponding result on \(\R^n\) due to P.-L. Lions \cite{Lions2a}*{Lemma I.1}.
P.-L.\thinspace Lions, F.\thinspace Pacella and M. Tricarico \cite{LionsPacellaTricarico1988}*{lemma 2.2} have adapted it to functions vanishing on a part of the boundary.

\begin{proof}
We follow the proof of P.-L. Lions \cite{Lions2a}*{lemma 1.1}. 
First assume that \(u=0\). Then, using Lemma~\ref{lemmaInequality} and Rellich's compactness theorem, one shows that for every \(a \in \Bar{\Omega}\) and \(\varepsilon > 0\), there exists \(\delta > 0\) such that if \(\varphi \in C(\Bar{\Omega})\), \(\varphi \ge 0\) and \(\supp \varphi \subset B(a, \delta)\), then 
\[
  \Bigl(\int_{\Bar{\Omega}} \varphi^{\frac{n}{n-1}} \mu \Bigr)^{1-\frac{1}{n}} \le \int_{\Bar{\Omega}} \varphi \nu.
\]
One deduces then the conclusion when \(u=0\) by the argument of \cite{Lions2a}*{lemma 1.2}.

The case \(u \ne 0\) follows then by standard arguments.
\end{proof}

\begin{proof}[Proof of proposition~\ref{propositionSufficientDomain}]
Let \((u_m)_{m \in \N}\) be a sequence in \(\mathrm{BV}(\Omega)\) such that 
\begin{align*}
  \int_{\Omega} \abs{Du_m} &\to c^q_\Omega,&
  \int_{\Omega} \abs{u_m}^{\frac{n}{n-1}} &= 1,&
  \int_{\Omega} \abs{u_m}^{q-1}u_m & = 0.
\end{align*}
Going if necessary to a subsequence, we can assume that the assumptions of lem\-ma~\ref{propositionCCLemma} are satisfied. Since \((u_m)_{m \in \N}\) converges weakly to \(u \in \mathrm{BV}(\Omega)\) and \(q < \frac{n}{n-1}\), by Rellich's compactness theorem, 
\[
  \int_{\Omega} \abs{u}^{q-1} u = \lim_{m \to \infty} \int_{\Omega} \abs{u_m}^{q-1}u_m = 0.
\]
Assume by contradiction that 
\[
  \int_{\Omega} \abs{u}^\frac{n}{n-1} < 1.
\]
In view of proposition~\ref{propositionCCLemma}, we have
\[
  \lim_{m \to \infty} \int_{\Omega} \abs{u_m}^\frac{n}{n-1}=\int_{\Omega} \nu = \int_{\Omega} \abs{u}^\frac{n}{n-1} + \sum_{j \in J} \nu_j,
\]
and thus \(J \ne \emptyset\). 
On the other hand,
\[
\begin{split}
 c^q_\Omega=\lim_{m \to \infty} \int_{\Omega} \abs{Du_m}=\int_{\Omega} \mu 
&\ge \int_{\Omega} \abs{Du} + \sum_{j \in J} \frac{c^*_n}{2^\frac{1}{n}} \nu_j^{1-\frac{1}{n}}\\
&> c^q_\Omega \Bigl(\int_{\Omega} \abs{u}^\frac{n}{n-1}\Bigr)^{1-\frac{1}{n}}+c^q_\Omega \sum_{j \in J} \nu_j^{1-\frac{1}{n}}\\
&\ge c^q_\Omega,
\end{split}
\]
which is a contradiction.
\end{proof}

\subsection{Upper estimate on the sharp constant}
We shall now prove that the condition of proposition~\ref{propositionSufficientDomain} is indeed satisfied.

\begin{proposition}
\label{propositionUpperEstimateDomain}
Let \(\Omega \subset \R^n\) be a bounded domain with \(\mathcal{C}^2\) boundary. If \(q < \frac{n^2}{n^2-1}\), then 
\[
  c^q_\Omega < \frac{c^*_n}{2^\frac{1}{n}}.
\]
\end{proposition}

\begin{proof}
Since \(\Omega\) is bounded, there exists \(a, b \in \partial \Omega\) such that \(\abs{a-b}=\sup \{\abs{x-y} \st x,\, y \in \partial \Omega\}\). Since \(\partial \Omega\) is of class \(\mathcal{C}^2\), its mean curvature \(H_a\) at \(a\) satisfies \(H_a \ge \frac{1}{\abs{a-b}} > 0\). For \(\varepsilon > 0\) such that \(\Omega \setminus \overline{B(a, \varepsilon)} \ne \emptyset\), consider the function \(u_\varepsilon \colon \Omega \to \R\) defined by
\[
  u_\varepsilon = \chi_{\Omega \cap B(a, \varepsilon)}-\beta_\varepsilon \chi_{\Omega \setminus B(a, \varepsilon)},
\]
where
\[
 \beta_\varepsilon =\Bigl(\frac{\muleb{n}(\Omega)}{\muleb{n}\bigl(\Omega \cap B(a, \varepsilon)\bigr)}-1 \Bigr)^{-\frac{1}{q}}.
\]
The quantity \(\mathcal{L}^n\bigl(\Omega \cap B(a, \varepsilon)\bigr)\) can be expanded in terms of the mean curvature \cite{HulinTroyanov2003}*{equation (1)} as
\[
 \muleb{n}\bigl(\Omega \cap B(a, \varepsilon)\bigr)= \frac{\pi^{\frac{n}{2}}\varepsilon^n}{2\Gamma(\frac{n}{2}+1)}\Bigl(1-\frac{n H_a\varepsilon}{(n+1) B(\tfrac{1}{2}, \tfrac{n-1}{2})} +o(\varepsilon)\Bigr),
\]
where \(B\) denotes Euler's beta function. 
In particular, one has
\[
 \beta_\varepsilon= \Bigl(\frac{\pi^{\frac{n}{2}}}{2\Gamma(\frac{n}{2}+1)\muleb{n}(\Omega)}\Bigr)^\frac{1}{q} \varepsilon^{\frac{n}{q}} \bigl(1+o(1)\bigr). 
\]

Since \(q < \frac{n^2}{n^2-1}\), we have \(\beta_\varepsilon^\frac{n}{n-1}=o(\varepsilon^{n+1})\) and therefore
\[
\begin{split}
   \int_{\Omega} \abs{u_\varepsilon}^{\frac{n}{n-1}}&=\muleb{n}\bigl(\Omega \cap B(a, \varepsilon)\bigr) - \beta_\varepsilon^{\frac{n}{n-1}} \muleb{n}\bigl(\Omega \setminus B(a, \varepsilon)\bigr) \\
&=\frac{\pi^{\frac{n}{2}}\varepsilon^n}{2\Gamma(\frac{n}{2}+1)}\Bigl(1-\frac{n}{n+1}\frac{H_a\varepsilon}{B(\tfrac{1}{2}, \tfrac{n-1}{2})}+o(\varepsilon)\Bigr).
\end{split}
\]
Similarly, one computes
\[
\begin{split}
 \int_{\Omega} \abs{Du_\varepsilon}&=(1+\beta_\varepsilon)\int_{\Omega} \abs{D\chi_{\Omega \cap B(a, \varepsilon)}} \\&=\varepsilon^{n-1}\frac{\pi^{\frac{n}{2}}n}{2\Gamma(\frac{n}{2}+1)}\Bigl(1-\frac{H_a \varepsilon}{B(\tfrac{1}{2}, \tfrac{n-1}{2} )}+o(\varepsilon)\Bigr).
\end{split}
\]
One has finally, since \(\beta_\varepsilon=o(\varepsilon)\).
\[
 \frac{\int_{\Omega} \abs{Du_\varepsilon}}{\Bigl(\int_{\Omega} \abs{u_\varepsilon}^{\frac{n}{n-1}}\Bigr)^{1-\frac{1}{n}}} 
=\frac{\pi^\frac{1}{2}n}{\bigl(2\Gamma(\frac{n}{2}+1)\bigr)^\frac{1}{n}}\Bigl(1- \frac{2 H_a \varepsilon}{(n+1)B(\frac{1}{2}, \frac{n-1}{2})} +o(\varepsilon)\Bigr),
\]
it follows then that for \(\varepsilon > 0\) sufficiently small, 
\[
 \frac{\int_{\Omega} \abs{Du_\varepsilon}}{\Bigl(\int_{\Omega} \abs{u_\varepsilon}^{\frac{n}{n-1}}\Bigr)^{1-\frac{1}{n}}} < c^*_n,
\]
which is the desired conclusion.
\end{proof}

In the previous proof, the existence of a point of the boundary with positive mean curvature is crucial. We would like to point out that in the problem of optimal functions for Sobolev--Hardy inequalities with a point singularity, one needs the boundary to have \emph{negative} mean curvature at that point of the boundary \cite{GhoussoubRobert}.

\begin{proposition}
\label{propositionMonotonicityDomains}
Let \(n \ge 2\) and \(\Omega \subset \R^n\) be a bounded domain with \(\mathcal{C}^1\) boundary.
For every \(q \in (0, \frac{n}{n-1})\), 
\[
  c^q_{\Omega} \le c^{\frac{1}{n-1}}_{\Omega}. 
\]
\end{proposition}
\begin{proof}
Let \(u \in \mathrm{BV}(\Omega)\) be such that \(\int_{\Omega} \abs{u}^{\frac{1}{n-1}-1}u=0\). Consider \(\lambda \in \R\) such that \(\int_{\Omega} \abs{u-\lambda}^{q-1}(u-\lambda)=0\). One has then 
\[
  c^{q}_\Omega \Bigl(\int_{\Omega} \abs{u}^{\frac{n}{n-1}}\Bigr)^{1-\frac{1}{n}}\le 
c^{q}_\Omega \Bigl(\int_{\Omega} \abs{u-\lambda}^{\frac{n}{n-1}}\Bigr)^{1-\frac{1}{n}}
 \le \int_{\Omega} \abs{D(u-\lambda)}=\int_{\Omega} \abs{Du},
\]
and therefore \(c^q_\Omega \le c^\frac{1}{n-1}_\Omega\).	
\end{proof}

By combining proposition~\ref{propositionUpperEstimateDomain} together with proposition~\ref{propositionMonotonicityDomains} we obtain

\begin{proposition}
\label{propositionUpperEstimateAllDomains}
Let \(n \ge 2\) and \(\Omega \subset \R^n\) be a bounded domain with \(\mathcal{C}^2\) boundary.
For every \(q \in (0, \frac{n}{n-1})\), 
\[
  c^q_\Omega < \frac{c^*_n}{2^\frac{1}{n}}. 
\]
\end{proposition}
\begin{proof}
Since \(n \ge 2\), one has
\[
  \frac{1}{n-1}  < \frac{n^2}{n^2-1}.
\]
Therefore, by proposition~\ref{propositionUpperEstimateDomain}, \(c^{\frac{1}{n-1}}_\Omega < \frac{c^*_n}{2^\frac{1}{n}}\). Hence, by proposition~\ref{propositionMonotonicityDomains}, 
\[
  c^q_\Omega \le c^\frac{1}{n-1}_\Omega < \frac{c^*_n}{2^{\frac{1}{n}}}.\qedhere
\]
\end{proof}

We are now in position to prove theorem~\ref{theoremDomainGeneral}, which contains theorems~\ref{theoremDomainAverage} and \ref{theoremDomainMinimal} as particular cases.

\begin{proof}[Proof of theorem~\ref{theoremDomainGeneral}]
By proposition~\ref{propositionUpperEstimateAllDomains}, \(c^q_\Omega < \frac{c^*_n}{2^\frac{1}{n}}\). Proposition~\ref{propositionSufficientDomain} is thus applicable and yields the conclusion.
\end{proof}

\section{Extremal functions on manifolds}

\subsection{Existence by concentration compactness}

The counterpart of proposition~\ref{propositionSufficientDomain} on manifolds is given by

\begin{proposition}
\label{propositionSufficientManifold}
Let \(n \ge 2\) and let \(M\) be an \(n\)--dimensional compact Riemannian \(\mathcal{C}^1\) manifold. If
\[
  c^q_M < c^*_n,
\]
then \(c^q_M\) is achieved.
\end{proposition}

\subsection{Upper estimate on the sharp constant}
We now turn on to estimates on the sharp constant,

\begin{proposition}
\label{propositionUpperEstimateManifold}
Let \(n \ge 3\) and \(M\) be an \(n\)--dimensional compact Riemannian \(\mathcal{C}^2\) manifold. If there exists \(a \in M\) such that the scalar curvature \(S_a\) at \(a\) is positive, then for every \(q \in (0, \frac{n^2}{n^2+n-2})\), 
\[
  c^q_M < c^*_n. 
\]
\end{proposition}

\begin{proof}
For \(\varepsilon > 0\) such that \(M \setminus \overline{B(a, \varepsilon)} \ne \emptyset\), where \(B(a,\varepsilon)\) is a geodesic ball of radius \(\varepsilon\) centered at \(a\), consider the function \(u_\varepsilon \colon M \to \R\) defined by
\[
  u_\varepsilon = \chi_{B(a, \varepsilon)}-\beta_\varepsilon \chi_{M \setminus B(a, \varepsilon)},
\]
where
\[
 \beta_\varepsilon =\Bigl(\frac{\muhaus{n}(M)}{\muhaus{n}\bigl(B(a, \varepsilon)\bigr)}-1 \Bigr)^{-\frac{1}{q}}.
\]
The measure of the ball can be extended as 
\[
 \muhaus{n}\bigl(B(a, \varepsilon)\bigr)=\frac{\pi^{\frac{n}{2}}\varepsilon^n }{\Gamma(\frac{n}{2}+1)}\Bigl(1-  \frac{S_a\varepsilon^2}{6(n+2)}+o(\varepsilon^2)\Bigr)
\]
(see for example \cite{Gray1973}*{Theorem 3.1}).
In particular, one has
\[
 \beta_\varepsilon= \Bigl(\frac{\pi^{\frac{n}{2}}}{\Gamma(\frac{n}{2}+1)\muhaus{n}(M)}\Bigr)^\frac{1}{q} \varepsilon^{\frac{n}{q}} \bigl(1+o(1)\bigr). 
\]

Since \(q < \frac{n^2}{n^2+n-2}\), we have \(\beta_\varepsilon^\frac{n}{n-1}=o(\varepsilon^{n+2})\) and therefore
\begin{equation}
\label{eqAsymptotuen}\int_{M} \abs{u_\varepsilon}^{\frac{n}{n-1}}=\frac{\pi^{\frac{n}{2}}\varepsilon^n}{\Gamma(\frac{n}{2}+1)}\Bigl(1-\frac{S_a\varepsilon^2}{6(n+2)}+o(\varepsilon^2)\Bigr).
\end{equation}
One also computes
\[
\begin{split}
 \int_{M} \abs{Du_\varepsilon}&=(1+\beta_\varepsilon)\int_{M} \abs{D\chi_{B(a, \varepsilon)}} \\
  &=\frac{n\pi^{\frac{n}{2}}\varepsilon^{n-1}}{\Gamma(\frac{n}{2}+1)}\Bigl(1-\frac{S_a \varepsilon^2}{6n}+o(\varepsilon^2)\Bigr).
\end{split}
\]
The combination of the previous developments gives
\begin{equation} 
 \label{eqAsymptotQuotientManifold}\frac{\displaystyle\int_{M} \abs{Du_\varepsilon}}{\displaystyle\Bigl(\int_{M} \abs{u_\varepsilon}^{\frac{n}{n-1}}\Bigr)^{1-\frac{1}{n}}} 
=\frac{n \pi^\frac{1}{2} }{\Gamma(\frac{n}{2}+1)^\frac{1}{n}} \Bigl(1- \frac{S_a \varepsilon^2}{2n(n+2)} +o(\varepsilon^2)\Bigr).
\end{equation}
It follows then that for \(\varepsilon > 0\) sufficiently small, 
\[
 \frac{\displaystyle\int_{M} \abs{Du_\varepsilon}}{\biggl(\displaystyle\int_{M} \abs{u_\varepsilon}^{\frac{n}{n-1}}\biggr)^{1-\frac{1}{n}}} < c^*_n,
\]
which is the desired conclusion.
\end{proof}

The counterpart of proposition~\ref{propositionMonotonicityDomains} can be obtained straightforwardly

\begin{proposition}
\label{propositionMonotonicityManifolds}
Let \(n \ge 2\) and \(M \subset \R^n\) be an \(n\)--dimensional compact Riemannian manifold.
For every \(q \in (0, \frac{n}{n-1})\), 
\[
  c^q_M \le c^{\frac{1}{n-1}}_M. 
\]
\end{proposition}

This allows us to obtain a counterpart of proposition~\ref{propositionUpperEstimateAllDomains}

\begin{proposition}
\label{propositionUpperEstimateAllnManifolds}
Let \(n \ge 3\) and \(M\) be an \(n\)--dimensional compact Riemannian \(\mathcal{C}^2\) manifold. If there exists \(a \in M\) such that the scalar curvature \(S_a\) at \(a\) is positive, then for every \(q \in (0, \frac{n}{n-1})\), 
\[
  c^q_M < c^*_n. 
\]
\end{proposition}
\begin{proof}
One checks that if \(n \ge 3\), \(\frac{1}{n-1} < \frac{n^2}{n^2+n-2}\). One can then apply proposition~\ref{propositionUpperEstimateManifold} with \(q=\frac{1}{n-1}\) and then conclude with proposition~\ref{propositionMonotonicityManifolds}.
\end{proof}

This allows us to prove theorem~\ref{theoremnManifolds} on manifolds

\begin{proof}[Proof of theorem~\ref{theoremnManifolds}]
Since \(n \ge 3\), this follows from proposition~\ref{propositionUpperEstimateAllnManifolds} and proposition~\ref{propositionSufficientManifold}.
\end{proof}

We can also prove theorem~\ref{theorem2ManifoldsSubcritical} on surfaces.

\begin{proof}[Proof of theorem~\ref{theorem2ManifoldsSubcritical}]
This follows from proposition~\ref{propositionUpperEstimateManifold} and proposition~\ref{propositionSufficientManifold}.
\end{proof}

\subsection{Refined upper estimates}
We now give a condition on the scalar curvature that gives a strict inequality in the critical case \(q=\frac{n^2}{n^2+n-2}\). Although we only need the result for \(n=2\), we state is for all dimensions.

\begin{proposition}
\label{propositionUpperEstimateManifoldCritical}
Let \(M\) be an \(n\)--dimensional compact Riemannian \(\mathcal{C}^2\) manifold. If there exists \(a \in M\) such that the scalar curvature \(S_a\) at \(a\) satisfies
\begin{equation}
\label{eqCriticalScalarCurvature}
 S_a > \frac{2(n+2)}{n-1}\Bigl(\frac{\pi^{\frac{n}{2}}\Gamma(\frac{n}{2}+1)}{\muhaus{n}(M)}\Bigr)^\frac{2}{n},
\end{equation}
then for \(q=\frac{n^2}{n^2+n-2}\),
\[
  c^q_M <  c^*_n.
\]
\end{proposition}

\begin{proof}
Since \(q=\frac{n^2}{n^2+n-2}\), the computations of the proof of proposition~\ref{propositionUpperEstimateManifold} give instead of \eqref{eqAsymptotuen}
\[
   \int_{M} \abs{u_\varepsilon}^{\frac{n}{n-1}}=\frac{\pi^{\frac{n}{2}}\varepsilon^n}{\Gamma(\frac{n}{2}+1)}\Bigl(1+\Bigl(\frac{\Gamma(\frac{n}{2}+1)}{\pi^{\frac{n}{2}}\muhaus{n}(M)}\Bigr)^\frac{2}{n}\varepsilon^2-\frac{S_a}{6(n+2)}\varepsilon^2+o(\varepsilon^2)\Bigr),
\]
and then, instead of~\eqref{eqAsymptotQuotientManifold},
\begin{multline*}
  \frac{\displaystyle\int_{M} \abs{Du_\varepsilon}}{\Bigl(\displaystyle\int_{\Omega} \abs{u_\varepsilon}^{\frac{n}{n-1}}\Bigr)^{1-\frac{1}{n}}} \\
=\frac{n\pi^\frac{1}{2}}{\Gamma(\frac{n}{2}+1)^\frac{1}{n}}\Bigl(1+\frac{n-1}{n}\Bigl(\frac{\Gamma(\frac{n}{2}+1)}{\pi^{\frac{n}{2}}\muhaus{n}(M)}\Bigr)^\frac{2}{n}\varepsilon^2- \frac{S_a}{2n(n+2)}\varepsilon^2 +o(\varepsilon)\Bigr),
\end{multline*}
and one checks that in view of \eqref{eqCriticalScalarCurvature}, the inequality is satisfied by taking some small \(\varepsilon > 0\).
\end{proof}

\begin{proof}[Proof of theorem~\ref{theorem2ManifoldsCritical}]
One first notes that proposition~\ref{propositionUpperEstimateManifoldCritical} is applicable with \(q=1\). Therefore by proposition~\ref{propositionMonotonicityManifolds}, for every \(q \in (0, 2)\), \(c^q_M < c^*_n\). The conclusion comes from proposition~\ref{propositionSufficientManifold}.
\end{proof}

\begin{proof}[Proof of theorem~\ref{theoremEulerPoincare}]
Since \(M\) does not have constant scalar curvature, there exists \(a \in M\) such that 
\[
S_a > \frac{1}{\muhaus{2}(M)} \int_M S.
\]
Since by the  Gauss--Bonnet formula
\[
  \int_M S = 4\pi \chi(M),
\]
we have
\[
  S_a > \frac{8\pi}{\muhaus{2}(M)}.
\]
The conclusion is then given by theorem~\ref{theorem2ManifoldsCritical}.
\end{proof}

\subsection{The case of the sphere}

\begin{proof}[Proof of theorem~\ref{theoremSphere}]
By proposition~\ref{propositionSufficientManifold}, we can assume that \(c^q_{\mathbf{S}^2} \ge c^*_2\). Let \(a \in \mathbf{S}^2\) and consider the function \(u \colon \mathbf{S}^2 \to \R\) defined by 
\[
  u= \chi_{B(a, \frac{\pi}{2})} - \chi_{\mathbf{S}^2 \setminus B(a, \frac{\pi}{2})},
\]
i.e.\ the difference between characteristic functions of opposite hemispheres.
One checks that 
\[
  \int_{\mathbf{S}^2} \abs{u}^{q-2} u =0 
\]
and
\[
 \frac{\displaystyle\int_{\mathbf{S}^2} \abs{Du}}{\Bigl( \displaystyle\int_{\mathbf{S}^2} \abs{u}^2 \Bigr)^\frac{1}{2}}=2 \sqrt{\pi}=c^*_2.
\]
Since we have assumed that \(c^q_{\mathbf{S}^2} \ge c^*_2\), this proves that the \(c^q_M\) is achieved.
\end{proof}

\end{document}